\documentclass[11pt]{amsart}

\renewcommand{\bar}{\overline}
\newcommand{\eps}{\varepsilon}
\newcommand{\pa}{\partial}
\renewcommand{\phi}{\varphi}

\newcounter{hours}\newcounter{minutes}

\newcommand{\M}{{\mathcal M}}

\newcommand{\ka}{K\"ahler }

\font\strange=msbm10

\newcommand{\C}{{{\mathchoice  {\hbox{$\textstyle{\text{\strange C}}$}}
{\hbox{$\textstyle{\text{\strange C}}$}}
{\hbox{$\scriptstyle{\text{\strange C}}$}}
{\hbox{$\scriptscriptstyle{\text{\strange C}}$}}}}}

\newcommand{\R}{{{\mathchoice  {\hbox{$\textstyle{\text{\strange R}}$}}
{\hbox{$\textstyle{\text{\strange R}}$}}
{\hbox{$\scriptstyle  N\kern-0.3em  R$}}  
{\hbox{$\scriptscriptstyle  R\kern-0.2em  R$}}}}}

\newcommand{\Z}{{{\mathchoice  {\hbox{$\textstyle{\text{\strange Z}}$}}
{\hbox{$\textstyle{\text{\strange Z}}$}}
{\hbox{$\scriptstyle  Z\kern-0.3em  Z$}}
{\hbox{$\scriptscriptstyle  Z\kern-0.2em  Z$}}}}}

\newcommand{\N}{{{\mathchoice  {\hbox{$\textstyle{\text{\strange N}}$}}
{\hbox{$\textstyle{\text{\strange N}}$}}
{\hbox{$\scriptstyle  N\kern-0.3em  N$}}
{\hbox{$\scriptscriptstyle  N\kern-0.2em  N$}}}}}

\newcommand{\bb}{{\frac{\sqrt{-1}}{2\pi}}}

\usepackage{amsmath,amsthm,amscd}
\usepackage{calc} 

\title 
[]{$K$ Energy and $K$ stability on Hypersurfaces }
\author{Zhiqin Lu}
\date{May 12, 2001}
\subjclass{Primary: 53A30; Secondary: 32C16}

\keywords{Mabuchi energy, Futaki invariants, and K-stability}
\address[Zhiqin Lu]
{Department of Mathematics\\
University of California, Irvine\\ 
Irvine, CA 92697}
\email[Zhiqin Lu]{zlu@math.uci.edu}
\thanks{Research supported by NSF grant DMS 9971506}

\newtheorem{theorem}{Theorem}[section]

\newtheorem{lemma}{Lemma}[section]
\newtheorem{cor}{Corollary}[section]
\newtheorem{prop}{Proposition}[section]

\newtheorem{definition}{Definition}[section]

\theoremstyle{remark}
\newtheorem{rem}{Remark}[section]

\begin{document}
\maketitle


\numberwithin{equation}{section}

\tableofcontents

\section{Introduction}
In this paper, we study the limiting property of the $K$ energy on
compact \ka hypersurface of $CP^n$ with positive first Chern class.

For a compact \ka manifold with positive first Chern class, one of the
most important problems is the existence of K\"ahler-Einstein
metric. If the manifold is a complex surface, the problem was
solved in~\cite{T6}. In higher dimensions, the existence of 
K\"ahler-Einstein metrics is 
related to certain geometric stability(cf.~\cite{T3}).

The notation $K$-stability was introduced in~\cite{T3} as
a necessary condition to the existence of K\"ahler-Einstein
metric. It is defined as follows:

Let $M$ be a compact \ka manifold in $CP^n$ such that there is
a constant $\alpha>0$ with $\alpha\omega_{FS}\in c_1(X)$. $M$ has this 
property if it is Fano and if the embedding is given by  the
anticanonical bundle. 
Let $\sigma(t)$ be a one parameter family of automorphisms
of $CP^n$. We write
\[
\sigma(t) [Z_0,\cdots,Z_n]=[t^{\lambda_0}Z_0,\cdots,
t^{\lambda_n}Z_n]
\]
for integers $\lambda_0,\cdots,\lambda_n$ with
$\sum\lambda_i=0$.
 Then we can define a 
family of metrics $\omega_t=\sigma(t)^*\omega_{FS}$ on $M$ such
that $\alpha\omega_t\in c_1(M)$.  Let $\M(\omega,\omega_t)$ be the
$K$ energy with respect to the metric $\alpha\omega$ and
$\alpha\omega_t$ (for the definition of the $K$ energy
defined by $K$, see next
section). It is known that
\begin{equation}\label{1-1}
\underset{t\rightarrow 0}{\rm lim}
\,t\frac{d}{dt}\M(\omega,\omega_t)=A
\end{equation}
exists~\cite{T3}. 
If $\M(\omega,\omega_t)$ has a lower bound,
then $A\geq 0$. Since the one parameter family of automorphisms
$\sigma(t)$
is generated by the holomorphic vector field $X=\sum\lambda_iZ_i
\frac{\pa}{\pa Z_i}$, we come up with
the following definition~\cite{T3}: 

\begin{definition}\label{d-11}
We say that $M$ is $K$ stable if for any holomorphic vector field
$X$ on $CP^n$ with $\lambda_0,\cdots,\lambda_n$ integers, 
\[
 \underset{t\rightarrow 0}{\rm lim}
\,t\frac{d}{dt}\M(\omega,\omega_t)> 0.
\]
If the above quantity is nonnegative for all vectors
$X$ on $CP^n$, we say $M$ is $K$
semistable.
\end{definition}

The general setting which relates the $K$ energy and the 
Futaki invariant is as follows: Let $M$ be a hypersurface
of $CP^n$. Let $X$ be the vector field of $CP^n$ in
Definition~\ref{d-11}. Suppose $M$ is defined by a polynomial
$F=0$ and
let $F_t=\sigma(-t)^*F$. The degeneration of $M$ by $X$ is defined
as the hypersurface in $\C\times CP^n$  by $G(t,\cdots)
=F_t(\cdots)=0$. The center fiber of the degeneration is defined
as the intersection of the degeneration with the set 
$\{0\}\times CP^n$, excluding the factor $t^\alpha=0$.

\begin{rem}
Definition~\ref{d-11} is a little bit more general than that
in~\cite{T3}. In fact, in~\cite{DT} or ~\cite{T3}, the quantity $A$ is
represented as the (real part) of the (generalized) Futaki invariant of
the center fiber if  the center fiber is a normal variety. However, the
exact same  proof can go through if we assume that the center fiber does
not have multiplicity greater than  1 (that implies, one can define the
``{\sl Futaki}'' invariant the same as the usual one for algebraic
cycles with multiplicity 1). 
\end{rem}

\begin{rem}
For our application, we only need the notion of $K$ semistability,
since
our first result only works on a dense subset of  all vector
fields on $CP^n$. Thus for the sake of simplicity, in this paper,
 we will use the terminology $K$ stability for both  $K$ stability
and
$K$ semistability. On the other hand, when we consider the 
limiting property of the $K$ energy, it doesn't make much
difference assuing $t$ is real or complex. Thus in the rest 
of this paper, we always assume that $t$ is a real number.
\end{rem}

The motivation of our work is to find an effective way to verify
the $K$ stability. In general, this is a harder problem than the
problem of finding an effective way to compute the Futaki invariant,
because the $K$ energy is the nonlinear version of the 
``Futaki'' invariant(see~\cite{Ma2}).
By the work of \cite{DT} or \cite{T3}, if the center fiber is
normal, the quantity $A$ is the real part of the corresponding Futaki
invariant. In this paper, we consider the case where $M$ is a 
hypersurface in $CP^n$ of degree less than or equal to $n$. Then
$M$ is a Fano variety and  one can compute the Futaki
invariant  by a very simple formula in ~\cite{Lu9} (see
also~\cite{Yo} by a completely different method). 

The technical difficulty in the proof is that the degeneration
of a hypersurface under a one parameter subgroup is ``generically''
an algebraic cycle of multiplication greater than 1. 
If that is the case, we would not be able to generalize
the argument in~\cite{DT} directly. In fact, our result
shows that the limit may {\sl not} depend on the center fiber
alone. This is on one side unexpected by the work of ~\cite{DT}.
On the other side, one realizes that if the center fiber
is a union of hyperplanes, it contains too little information of
the degeneration so that extra piece of information 
from the degeneration is needed. 

In this paper, we
overcome the above difficulty
in the case that the center fiber is of multiplicity great than one.
We first represent the $K$ energy into an
explicitly formula(Theorem~\ref{thm21}). Then we carefully analyze the
integrand in the formula  by using some analytic techniques and 
a recently result of Phong and Sturm~\cite{PS} to get the conclusion.

Our result can be generalized to case of complete intersections
and even general projective manifold. The results will appear 
in a subsequent paper.

Before stating the main result, we setup some notations:
let $M$ be
defined by the zeros of the polynomial
\begin{equation}\label{pp-1}
F=\sum_{i=0}^p a_iZ_0^{\alpha_0^i}\cdots Z_n^{\alpha_n^i}
\end{equation}
of degree $d$.
Let $(\lambda_0,\cdots,\lambda_n)$ be the rational numbers 
satisfying $\sum\lambda_i=0$. Let
\begin{equation}\label{pp-2}
\lambda=\underset{0\leq i\leq p}{{\rm Max}}
(\sum_{k=0}^n\lambda_k\alpha_k^i).
\end{equation}
Let
\begin{equation}\label{pp-3}
\phi(x_0,\cdots,x_n)=\underset{0\leq i\leq p}
{\rm Min} (-\sum_{k=0}^n\lambda_k
\alpha_k^i
+\sum_{k=0}^n\alpha_k^i
x_k),
\end{equation}
and let
\begin{equation}\label{pp-4}
\phi_i(x)=\phi(0,\cdots,\underset{i}{x},\cdots,0).
\end{equation}
Then we have the following

\begin{theorem}\label{main}
For ``generic'' (See section 3 for details)
$(\lambda_0,\cdots,
\lambda_n)$, we have
\begin{align*}
&\qquad\underset{r\rightarrow 0}{\rm lim}\,t\frac{d}{dt}\M(t)\\&=
\frac{2}{d}\left(
-\frac{\lambda(d-1)(n+1)}{n}
+\sum_{i=0}^n\int_0^\infty\phi_i'(x)(
\phi_i'(x)-1)dx\right).
\end{align*}
\end{theorem}

Since for a K\"ahler-Einstein manifold, the $K$
energy has a lower bound, we have the following:

\begin{theorem}\label{mainc}
 If $M$ is a K\"ahler-Einstein 
hypersurface with positive first Chern class, then
 we have
\[
-\frac{\lambda(d-1)(n+1)}{n}
+\sum_{i=0}^n\int_0^\infty\phi_i'(x)(
\phi_i'(x)-1)dx\geq 0
\]
for any $\lambda_0\cdots,\lambda_n\in\R$ with $\sum\lambda_i=0$.
\end{theorem}

{\bf Proof of Theorem~\ref{mainc}.} The expression in the theorem
is continuous and homogeneous with respect to
$\lambda_0,\cdots,\lambda_n$. So by taking
the limit, we proved that the inequality is valid for any choice of 
$\lambda_0,\cdots,\lambda_n\in\R$. \qed

\begin{rem}
As proved in Theorem~\ref{normal}, one can define the generalized
Futaki invariant for the  degeneration as
$\underset{t\rightarrow 0}{\rm lim}\,t\M'(t)$. When the center fiber
is a normal variety, it is the generalized Futaki invariant
defined in~\cite{DT}  by the work of~\cite{Lu9}. However, Theorem
~\ref{mainc} is more powerful in the hypersurface case 
than that in~\cite{DT} because of the
flexibility of the choices of
$X$. 
\end{rem}

{\bf Acknowledgment.} The author thanks P. Li, D.H. Phong and G. Tian
for the encouragement during the preparation of this paper. Special
thanks to L. Katzarkov who helps the author clarify a lot
of concepts in algebraic geometry.

\section{An explicit formula for the $K$ energy}
In this section, we  give an explicit
formula for the $K$ energy of smooth hypersurfaces 
of $CP^n$. 

First, let's recall the definition of the $K$ energy
~\cite{Ma2}.
Let $M$ be a compact \ka manifold with positive first 
Chern class $c_1(M)$. Let $\omega_0,\omega_1
\in c_1(M)$ and let $\omega_1=\omega_0+\bb\pa\bar\pa\xi$
for a smooth function $\xi$. We
put
$\omega_s=\omega_0+s\bb\pa\bar\pa\xi$ and define
\begin{equation}\label{2-0} 
\M(\omega_0,\omega_1)
=-\frac{1}{V}\int_0^1\left(
\int_X\xi(R(\omega_s)-m)\omega_s^m\right)ds,
\end{equation}
where $R(\omega_s)$ is  the scalar curvature of the
metric, $m$ is the complex dimension of $M$  and $V$ is the
volume of $X$ with respect to $\omega_0$. The functional 
$\M$, which is called the $K$ energy by Mabuchi, has
the properties:
\begin{enumerate}
\item $\M(\omega_0,\omega_1)=-\M(\omega_1,\omega_0)$,
\item $\M(\omega_0,\omega_1)+\M(\omega_1,\omega_2)
=\M(\omega_0,\omega_2)$,
\end{enumerate}
where $\omega_0,\omega_1,\omega_2\in c_1(X)$.
 
From now on, let's assume that $\omega$ is the \ka form of
the Fubini-Study metric of
$CP^n$. Let $M$ be a hypersurface in $CP^n$ defined by the 
polynomial $F=0$ of degree $d$. 
Of course, we need $d\leq n$ to insure that $M$ is Fano.
Let
$\lambda_0,\cdots,\lambda_n$ be integers such that
$\sum_{i=0}^n\lambda_i=0$. 
Let $F_t$ be the polynomial
defined by
\[
F_t(Z_0,\cdots,Z_n)=
F(t^{-\lambda_0}Z_0,\cdots,t^{-\lambda_n}Z_n),
\]
and let $M_t$ be the hypersurface defined by the zero set
of $F_t$. Geometrically, $M_t$ is the image of $M$ 
under the
automorphism $\sigma(t)$ generated by the holomorphic vector
field
$X=\sum_{i=0}^n\lambda_iZ_i\frac{\pa}{\pa Z_i}$. The
automorphisms $\sigma(t)$ can be written as
$\sigma(t)([Z_0,\cdots,Z_n])=[
t^{\sigma_0}Z_0,\cdots,t^{\sigma_n}Z_n]$. Using these
automorphisms, one can define a family of \ka forms
$\omega_t=
\sigma(t)^*\omega$ on $M$. It is easy to see that both
$(n-d+1)\omega$ and $(n-d+1)\omega_t$ are \ka forms of $M$
in the cohomological class $c_1(M)$. Define
$\M(t)=\M((n-d+1)\omega,(n-d+1)\omega_t)$. It is a well
known result~\cite{Ma2} that if $M$ admits a
K\"ahler-Einstein metric, then
$\M(t)$ has a lower bound.

\begin{prop}\label{prop21} 
Using the notations as above, we have
\[
t\frac{d}{dt}\M(t)
=\frac{2(n-1)}{d}
\int_{M_t}({\rm
Ric}(\omega|_{M_t})-(n-d+1)\omega|_{M_t})\theta\omega^{n-2},
\]
where $\theta$ is defined as
\begin{equation}\label{3-1-1}
\theta=-\frac{\sum_{i=0}^n\lambda_i|Z_i|^2}
{\sum_{i=0}^n|Z_i|^2},
\end{equation}
and ${\rm Ric}(\omega|_{M_t})$ is the Ricci curvature
of $\omega|_{M_t}$.
\end{prop}
 
{\bf Proof.} It basically follows from the 
two properties
of the $K$ energy three lines under the
equation~\eqref{2-0}. See ~\cite[Lemma 2.1]{DT} for details.

\qed

The following lemma can be found in ~\cite{T8}, we include
the proof here for the sake of completeness. 
\begin{lemma}\label{lem3-1}
Let $M$ be the smooth hypersurface defined as the zero
of $\{F=0\}$. 
We use $\omega$ to denote the Fubini-Study metric on $CP^n$
as well as the K\"ahler form on $M$, which is the
restriction of $\omega$ on $M$.   
 Let
\begin{equation}\label{3-1}
\xi=\log\frac{|\nabla F|^2}{(\sum_{i=0}^n|Z_i|)^{2(d-1)}},
\end{equation}
where $[Z_0,\cdots, Z_n]$ is the homogeneous
coordinate in $CP^n$.
Then we have
\begin{equation}\label{result-3-1}
Ric (\omega)-(n-d+1)\omega=-\bb\pa\bar\pa\xi.
\end{equation}
\end{lemma}

{\bf Proof.} 
Without losing generality, 
we prove the above lemma on  the open set
$U_0=\{[Z_0,\cdots, Z_n]||Z_0|>\frac 12|Z_j|,
j=1,\cdots,n\}$ in
$CP^n$.  The local coordinate system on $U_0$ is
$(z_1,\cdots,z_n)$ where $z_i=Z_i/Z_0$ for $i=1,\cdots n$. Under this 
coordinate system, the Fubini-Study metric can be written as
\begin{equation}\label{3-2}
\omega=\bb g_{i\bar j}dz_i\wedge d\bar{z}_j
=\bb\sum_{i,j=1}^n(\frac{\delta_{ij}}{1+|z|^2}-\frac{z_j\bar{z}_i}
{(1+|z|^2)^2})dz_i\wedge d\bar{z}_j,
\end{equation}
where $|z|^2=\sum |z_i|^2$. Let's further assume that  in a small open
set $V$ of
$U_0$, from the equation $F=0$, we can solve $z_1$. Namely,
\begin{equation}\label{3-3}
z_1=z_1(z_2,\cdots,z_n)
\end{equation}
for a holomorphic function $z_1$.
Let the \ka form $\omega$ on $V$, under the local
coordinate system $(z_2,\cdots, z_n)$, be written as
\[
\omega=\bb\sum_{i,j=2}^n\tilde{g}_{i\bar j} dz_i\wedge
d\bar z_j,
\]
and let $a_i=\frac{\pa z_1}{\pa z_i}, i=2,\cdots ,n$.
Then by~\eqref{3-2} and~\eqref{3-3}, we have
\begin{align*}
 &\tilde g_{i\bar j}=\frac{\delta_{ij}}{1+|z|^2}-\frac{z_j\bar{z}_i}
{(1+|z|^2)^2}-
\frac{z_j\bar{z_1}a_i}{(1+|z|^2)^2}
-\frac{z_1\bar{z}_i\bar{a}_j}{(1+|z|^2)^2}\\&\quad
+\frac{a_i\bar{a}_j}{1+|z|^2}-\frac{|z_1|^2 a_i\bar{a}_j}
{(1+|z|^2)^2},
\end{align*}
for $i,j=2,\cdots,n$. We want to compute the determinant $\det
(\tilde g_{i\bar j})$. In order to do  this, we let
\[
K_{i\bar j}=\delta_{ij}+a_i\bar{a}_j-\frac{1}
{1+|z|^2}(\bar z_i+\bar z_1a_i)\overline{(\bar z_j+\bar z_1 a_j)}.
\]
Then 
\begin{equation}\label{3-4}
\tilde g_{i\bar j}=\frac{1}{1+|z|^2}K_{i\bar j},\quad i,j=2,\cdots n.
\end{equation}

Let
\begin{align*}
&A=(a_2,\cdots, a_n);\\
&B=(\bar z_2+\bar z_1 a_2,\cdots,\bar z_n+\bar z_1a_n).
\end{align*}
Then the matrix $K=(K_{i\bar j})$ can be represented by
\[
K=I+A^T\bar A-\frac{1}{1+|z|^2}B^T\bar B.
\]
A straightforward computation gives
\begin{align*}
&KA^T=(1+|a|^2)A^T-\frac{1}{1+|z|^2}(\bar BA^T)B^T;\\
& KB^T=(\bar AB^T)A^T+(1-\frac{|B|^2}{1+|z|^2})B^T.
\end{align*}
Thus the vector space spanned by the vectors $A$, $B$ is $K$-invariant.
Furthermore, on the complement of the vector space, $K$ is the
identity. So we have
\begin{align}\label{3-5}
\begin{split}
&
\det K=(1+|a|^2)(1-\frac{|B|^2}{1+|z|^2})+\frac{1}{1+|z|^2}
|\bar BA^T|^2\\
&=\frac{1}{1+|z|^2}(1+|a|^2+|\sum_{i=2}^na_iz_i-z_1|^2).
\end{split}
\end{align}
Let $f$ be the defining function of $M$ on $U_0$, i.e.
\[
f=F(1,\frac{Z_i}{Z_0},\cdots,\frac{Z_n}{Z_0})=\frac{F}{Z_0^d}.
\]
Then
\begin{equation}\label{3-7}
\frac{\pa z_1}{\pa z_k}=-\frac{\frac{\pa f}{\pa z_k}}
{\frac{\pa f}{\pa z_1}}=-\frac{F_k}
{F_1},\quad (k=2,\cdots,n)
\end{equation}
where we define $F_k=\frac{\pa F}{\pa Z_k}$ for 
$k=0,\cdots,n$.
Thus by the homogeneity of $F$, we have
\begin{align}\label{3-7-8}
\begin{split}
&\sum_{i=2}^na_iz_i-z_1=-\sum_{i=2}^n\frac{Z_i}{Z_0}
\frac{F_i}{F_1}
-\frac{Z_1}{Z_0}\\
&=-\frac{1}{Z_0F_1}
(\sum_{i=1}^nZ_iF_i)
=\frac{F_0}{F_1}
\end{split}
\end{align}
on $M$.
Using~\eqref{3-4} and ~\eqref{3-5}, we have
\begin{equation}\label{3-7-1}
\det \tilde g_{i\bar j}=\frac{1}{(1+|z|^2)^n}\frac{1}
{|F_1|^2}(\sum_{k=0}^n|F_k|^2).
\end{equation}
Then by ~\eqref{3-1}
\[
\det \tilde g_{i\bar j}=\frac{1}
{(1+|z|^2)^{n-d+1}}\cdot\frac{1}
{|\frac{\pa f}{\pa z_1}|^2}\cdot e^\xi.
\]
\eqref{result-3-1} follows from the formula of the Ricci
curvature and 
the above equation.

\qed

In order to represent the $K$ energy in terms of the
polynomial $F$, we need the following
 purely algebraic lemma:

\begin{lemma}\label{lem32}
With the same notations as above, let $\eta$ be a $(1,1)$ form
on $CP^n$. Let $\pi: C^{n+1}\rightarrow CP^n$ be the projection.
Let
\begin{equation}\label{3-7-3}
\pi^*\eta=\bb\sum_{i,j=0}^{n}\tilde a_{i\bar j}dZ_i\wedge d\bar Z_j.
\end{equation}
Then on $M$,
\begin{equation}\label{3-7-3-1}
\eta\wedge\omega^{n-2}=\frac{|Z|^2}
{n-1}
\left(
\sum_{i=0}^n\tilde a_{i\bar i}-\frac{\sum_{i,j=0}^n\tilde a_{i\bar
j}F_j\bar F_i} {|\nabla F|^2}\right)\omega^{n-1}
\end{equation}
for $|Z|^2=\sum_{i=0}^n|Z_i|^2$.
\end{lemma}

\begin{rem}
The righthanded side of ~\eqref{3-7-3-1} is well defined 
because $\tilde a_{i\bar j}$ for $i,j=0,\cdots,n$ are
homogeneous functions of order $(-2)$.
\end{rem}

{\bf Proof.}
As in the proof of the previous lemma, we can consider the problem
only on $U_0\cap\{\frac{\pa F}{\pa Z_1}\neq 0\}$, without losing
generality. Define $A_{i\bar j}$ on $CP^n$ as follows:
\begin{align}
\begin{split}\label{3-6}
&
\eta\wedge\omega^{n-2}
=(\bb)^{n-1}(-1)^{\frac 12(n-1)(n-2)}\\
&\cdot
\sum_{i,j=1}^n (-1)^{i+j}A_{i\bar j}
dz_1\wedge\cdots\wedge\hat{dz^i}\cdots \wedge dz^n
\wedge
d\bar z_1\wedge\cdots \wedge\hat{d\bar z_j}\cdots \wedge d\bar z_n
,\end{split}
\end{align} 
where ``$\hat{\quad}$'' means omit.
Define
\[
b=(1, -a_2,\cdots, -a_n)
=(1,-\frac{\pa z_1}{\pa z_2},\cdots,
-\frac{\pa z_1}{\pa
z_n}) =(1,\frac{F_2}{F_1},
\cdots,
\frac{F_n}{F_1}).
\]
Then
by~\eqref{3-6}, we have 
\begin{align}
\begin{split}\label{3-7-2}
&
\eta\wedge\omega^{n-2}
=(\bb)^{n-1}(-1)^{\frac 12(n-1)(n-2)}\\
&
\cdot\sum_{i,j=1}^n A_{i\bar j}b_i\bar b_j
dz_2\wedge\cdots\wedge dz^n
\wedge
d\bar z^2\wedge
\cdots \wedge d\bar z^n
\end{split}
\end{align}
on $M$. Thus in order to prove~\eqref{3-7-3-1},
we just need to compute $\sum A_{i\bar j}b_i\bar b_j$. To this end,
let 
\begin{equation}\label{3-7-4}
\eta=\bb\sum_{i,j=1}^n  a_{k\bar l} dz_k\wedge d\bar z_l,
\end{equation}
and fix $r$, $s$. By~\eqref{3-6}, we have
\begin{align}
\begin{split}\label{3-8}
&
\bb dz_r\wedge d\bar z_s
\wedge
\bb
a_{k\bar l}dz_k\wedge d\bar z_l\wedge\omega^{n-2}\\
&
=(\bb)^n(-1)^{\frac 12(n-1)(n-2)}(-1)^{n-1}
A_{r\bar s}dz_1\wedge\cdots\wedge d\bar z_n.
\end{split}
\end{align}
We also have the following algebraic fact:
\begin{align}
\begin{split}\label{3-9}
&
\bb dz_r\wedge d\bar z_s
\wedge
\bb
a_{k\bar l}dz_k\wedge d\bar z_l\wedge\omega^{n-2}\\
&
=\frac{1}{n(n-1)}\left(
\sum_{\alpha, \beta=1}^n
(g^{\alpha\bar\beta}a_{\alpha\bar\beta})
g^{r\bar s}-
\sum_{\alpha,\beta=1}^ng^{\alpha\bar s}g^{r\bar\beta}
a_{\alpha\bar\beta}\right)\omega^n.
\end{split}
\end{align}
By~\eqref{3-2}, we have
\begin{equation}\label{3-10}
\omega^n=\left(\bb\right)^nn!(-1)^{\frac 12n(n-1)}
\frac{1}{(1+|z|^2)^{n+1}}dz_1\wedge\cdots\wedge d\bar z_n.
\end{equation}
Comparing ~\eqref{3-8}, ~\eqref{3-9} and~\eqref{3-10}, we have
\begin{equation}\label{3-11}
A_{r\bar s}=\frac{(n-2)!}{(1+|z|^2)^{n+1}}
(\sum_{\alpha, \beta=1}^ng^{\alpha\bar\beta}a_{\alpha\bar\beta}
g^{r\bar s}-\sum_{\alpha,\beta=1}^ng^{\alpha\bar s}g^{r\bar\beta}
a_{\alpha\bar\beta}),
\end{equation}
for $r,s=1,\cdots, n$. By~\eqref{3-11}, we have
\begin{align}\label{3-11-1}
\begin{split}
&\sum_{i,j=1}^nA_{i\bar j}b_i\bar b_j
=
\frac{(n-2)!}{(1+|z|^2)^{n+1}}\\&\cdot
(\sum_{\alpha,\beta=1}^n
g^{\alpha\bar\beta}
a_{\alpha\bar\beta}\sum_{i,j=1}^ng^{i\bar j}b_i\bar b_j
-\sum_{i,j,\alpha,\beta=1}^n
g^{\alpha\bar j}g^{i\bar\beta}a_{\alpha\bar\beta}b_i\bar
b_j).
\end{split}
\end{align}

We need the following
\begin{lemma}\label{lem322}
Using the same notations as above, we have
\begin{align}
&\sum_{\alpha,\beta=1}^n g^{\alpha\bar\beta} a_{
\alpha\bar\beta}=|Z_0|^2(1+|z|^2)\sum_{i=0}^n\tilde a_{i\bar
i},\label{3-15-1}\\
&\sum_{i,j=1}^n g^{i\bar j}b_i\bar b_j
=
(1+|z|^2)\frac{|\nabla F|^2}{|F_1|^2},\label{3-15-2}\\
&\sum_{i,j,\alpha,\beta=1}^n g^{\alpha\bar j} g^{i\bar
\beta}
a_{\alpha\bar\beta}b_i\bar b_j=
|Z_0|^2(1+|z|^2)^2\frac{
\sum_{\alpha,\beta=0}^n\tilde
a_{\alpha\bar\beta}
\bar F_\alpha F_\beta}{|F_1|^2},\label{3-15-3}
\end{align}
where $\tilde a_{i\bar j}$ is defined in~\eqref{3-7-3}.
\end{lemma}

{\bf Proof.} Comparing ~\eqref{3-7-3} and~\eqref{3-7-4}, we
have
\begin{equation}\label{3-15}
\left\{
\begin{array}{ll}
 a_{k\bar l}=\tilde a_{k\bar l}\cdot |Z_0|^2,& k,l\neq 0;\\
\sum_{i=1}^n z_i a_{i\bar l}=-\tilde a_{0,\bar l}\cdot
|Z_0|^2, &l\neq 0;\\
\sum_{j=1}^n\bar z_j a_{k\bar j}=
-\tilde a_{k\bar 0}\cdot |Z_0|^2,& k\neq 0;\\
 \sum_{i,j=1}^n z_i\bar z_j a_{i\bar j}
=\tilde a_{0\bar 0}\cdot |Z_0|^2.
\end{array}
\right.
\end{equation}
Since
$g^{\alpha\bar\beta}=(1+|z|^2)(\delta_{\alpha\beta}+z_\alpha\bar
z_\beta)$, by~\eqref{3-15}
, we have
\[
\sum_{\alpha,\beta=1}^n g^{\alpha\bar\beta}a_{\alpha\bar\beta}
=(1+|z|^2)\sum_{\alpha,\beta=1}^n
(\delta_{\alpha\beta}+z_\alpha\bar z_\beta)a_{\alpha\bar\beta}
=|Z|^2\sum_{\alpha=0}^n\tilde a_{\alpha\bar\alpha}.
\]
This proves~\eqref{3-15-1}. By~\eqref{3-7-8}, we have
\[
\sum_{i=1}^n z_ib_i=-\frac{F_0}{F_1}
\]
on $M$.
Thus
~\eqref{3-15-2}
and~\eqref{3-15-3} follow from a straightforward computation
using the above equation.

\qed

{\bf Continuation of the Proof of Lemma~\ref{lem32}.}
By Lemma~\ref{lem322}, we have
\begin{align}\label{3-12}
\begin{split}
&\sum_{\alpha,\beta=1}^ng^{\alpha\bar\beta}
a_{\alpha\bar\beta}
\sum_{i,j=1}^ng^{i\bar
j}b_i\bar b_j -
\sum_{i,j,\alpha,\beta=1}^n
g^{\alpha\bar j}g^{i\bar\beta}a_{\alpha\bar\beta}b_i\bar
b_j\\&
=
|Z_0|^2(1+|z|^2)^2
\frac{|\nabla F|^2}{|F_1|^2}\left(
\tilde a_{i\bar i}-\frac{\tilde a_{i\bar j}F_j\bar F_i}
{|\nabla F|^2}\right).
\end{split}
\end{align}

By ~\eqref{3-7-1},
\begin{align}
\begin{split}\label{3-13}
&\omega^{n-1}
=
\left(\bb\right)^{n-1}(-1)^{\frac 12(n-1)(n-2)}
\frac{(n-1)!}{(1+|z|^2)^n}\frac{|\nabla F|^2}{|F_1|^2}
\\
&
\qquad\cdot
dz_2\wedge\cdots\wedge dz_n\wedge d\bar z_2\wedge\cdots d\bar
z_n.
\end{split}
\end{align}
~\eqref{3-7-3-1} follows from
~\eqref{3-7-2}, ~\eqref{3-11-1}, \eqref{3-12}
and~\eqref{3-13}.

\qed

\begin{lemma}\label{lem310}
Let $\xi$ be the function defined in ~\eqref{3-1}
and let $\theta$ be defined in ~\eqref{3-1-1}. 
Then we have
\begin{align}\label{3-15-4}
\begin{split}
&\qquad
\bb\pa\xi\wedge\bar\pa\theta\wedge\omega^{n-2}\\
&=
\frac{1}{n-1}
\left(-\sum_{k=0}^n\left(\frac{XF}{|\nabla
F|^2}\right)_k\bar F_k+\frac{\sum_{k=0}^n\lambda_k|F_k|^2}
{|\nabla F|^2}-(d-1)\theta\right)\omega^{n-1}.
\end{split}
\end{align}
Furthermore, we have
\begin{align}
\begin{split}\label{3-14}
&\qquad\bb
\int_M\pa\xi\wedge\bar\pa\theta\wedge\omega^{n-2}\\
&=-\frac{1}{n-1}\int_M
\sum_{k=0}^n\left(\frac{XF}{|\nabla F|^2}\right)_k
\bar F_k\omega^{n-1}
+\frac{n-d+1}{n-1}\int_M\theta\omega^{n-1}.
\end{split}
\end{align}
\end{lemma}

{\bf Proof.} Let $\bb\eta=\pa\xi\wedge\bar\pa\theta$
and let 
\[
\pi^*\eta=\bb\sum_{i,j=0}^n\tilde a_{i\bar j} dZ_i\wedge
d\bar Z_j.
\]
Then we have
\[
\tilde a_{i\bar j}=\frac{\pa\xi}{\pa Z_i}\cdot
\frac{\pa\theta}{\pa\bar Z_j}.
\]
A straightforward computation gives
\[
\sum_{i=0}^n\tilde a_{i\bar i}=\frac
{-\sum_{k=0}^n(XF)_k\bar F_k+\sum_{k=0}^n \lambda_k|F_k|^2}
{|Z|^2|\nabla F|^2}-(d-1)\frac{\theta}{|Z|^2},
\]
and
\[
\frac{\sum_{i,j=0}^n\tilde a_{i\bar j}F_j\bar F_i}
{|\nabla F|^2}
=-\frac
{XF\cdot\sum_{i,k=0}^n F_{ik}\bar F_i\bar F_k}
{|Z|^2|\nabla F|^4}
\]
on $M$.
Thus by Lemma~\ref{lem32},
we got~\eqref{3-15-4}.
Using Lemma~\ref{lem32} again by setting
$\bb\eta=\pa\bar\pa\theta$, we have
\begin{equation}\label{3-16}
\bb\pa\bar\pa\theta\wedge\omega^{n-1}=\frac{1}{n-1}
(-n\theta+\frac{
\sum_{k=0}^n\lambda_k|F_k|^2}{|\nabla F|^2})
\omega^{n-1}.
\end{equation}
\eqref{3-14}  follows from ~\eqref{3-15-4}, ~\eqref{3-16}
and the Stokes Theorem.

\qed

Although  not needed in this paper, we give a simple 
proof of 
the following formula for the Futaki invariant
in~\cite{Lu9} as an application of Lemma~\ref{lem3-1}, 
Lemma~\ref{lem32} and Lemma~\ref{lem310}.

\begin{cor}\label{cor21}
Let $M$ be a smooth hypersurface in $CP^n$ defined
by the homogeneous polynomial $F=0$ of degree $d$. Let
$X$ be a vector in $CP^n$ satisfying
\begin{equation}\label{3-17}
XF=\kappa F.
\end{equation}
The Futaki invariant is defined as
\[
{\mathcal F}(X)=-\int_M X(\xi)\omega^{n-1}.
\]
Then
\begin{equation}\label{f1}
{\mathcal F}(X)=-\frac{(n+1)(d-1)}{n}\kappa.
\end{equation}
\end{cor}

{\bf Proof.}
We have
\begin{equation}\label{d3-1}
i(X)\omega=-\bar\pa\theta.
\end{equation}
Since $X$ leave $M$ invariant, we have
\[
0=\int_M i(X)(\pa\xi\wedge\omega^{n-1})
=
\int_MX\xi\omega^{n-1}+(n-1)\int_M
\pa\xi\wedge\bar\pa\theta\wedge\omega^{n-2}.
\]
By
the above equation and
~\eqref{3-17},
we have
\[
{\mathcal
F}(X)=-\int_M\omega^{n-1}+(n-d+1)\int_M\theta\omega^{n-1}.
\]
By~\cite[Theorem 5.1]{Lu9}, we have
\[
\int_M\theta\omega^{n-1}=\frac{\kappa}{n}.
\]
\eqref{f1} follows from the above two equations.

\qed  

Finally, we have the following
\begin{theorem}\label{thm21}
The $K$ energy $\M(t)$ can be
represented as
\begin{align}\label{m-1}
\begin{split}
&\qquad\M(t)=
\frac 2d\int_0^t\left(
\int_{M_\tau}\frac 1\tau\left(
-\sum_{k=0}^n\left(\frac{XF_\tau}{|\nabla
F_\tau|^2}\right)_k
\bar{(F_\tau)_k}\omega^{n-1}\right.\right.\\
&\qquad \left.\left.+
(n-d+1)\int_{M_\tau}\theta\omega^{n-1}
\right)\right)d\tau,
\end{split}
\end{align}
where
\[
F_\tau(Z_0,\cdots,Z_n)=F(\tau^{-\lambda_0}Z_0,\cdots,
\tau^{-\lambda_n}Z_n),
\]
and $M_\tau$ is the zero set of $F_\tau=0$. In particular,
we have
\begin{align}\label{m-2}
\begin{split}
&\qquad t\frac{d}{dt}\M(t)=\frac{2}{d}
\left(
-\sum_{k=0}^n\left(\frac{XF_t}{|\nabla
F_t|^2}\right)_k
\bar{(F_t)_k}\omega^{n-1}+
(n-d+1)\int_{M_t}\theta\omega^{n-1}\right).
\end{split}
\end{align}
\end{theorem}

{\bf Proof.} The theorem follows from
Prop~\ref{prop21}, Lemma~\ref{lem3-1} and
Lemma~\ref{lem310}.

\qed

\section{The limit of the derivative of the $K$ energy}
In this section, we compute the limit 
$\underset{t\rightarrow 0}{\rm lim}\,t\M'(t)$ using 
Proposition~\ref{prop21}. First,
we need  some combinatoric preparations.

Let $(\delta_i,\sigma_i), i=0,\cdots,p$ be a sequence of pair of
nonnegative rational numbers. $\delta_0=0$. We assume
that the sequence is ``generic'' in the sense that
\begin{enumerate}
\item All $\delta_i, (i=0,\cdots, p)$ are distinct numbers(that implies $\delta_i>0,
i=1,\cdots,p$);
\item None of the three lines defined by $\psi_i(x)=\delta_i+\sigma_ix,
(i=0,\cdots,p)$ intersect
at the same point.
\end{enumerate}

Define $(i_k, r_k), (k=0,\cdots,m)$ inductively as follows: let 
$i_0=0, r_0=0$. If $(i_k,r_k)$ has been defined, then 
\begin{enumerate}
\item If for any $r>r_k$
\[
\delta_{i_k}+\sigma_{i_k}r<\delta_i+\sigma_ir\qquad (i\neq i_k),
\]
then let $m=k$ and stop;
\item If not, then define $i_{k+1}$ and $r_{k+1}>r_k$ such that
\begin{equation}\label{4-1}
\delta_{i_k}+\sigma_{i_k}r_{k+1}=\delta_{i_{k+1}}
+\sigma_{i_{k+1}}r_{k+1}\leq \delta_i+\sigma_ir_{k+1},
\end{equation}
where $i=1,\cdots,p$. Since $(\delta_i,\sigma_i), i=0,\cdots,p$ are
 ``generic", the
choice of $(i_k, r_k)$ is unique for
 $(k=0,\cdots,m)$ .
\end{enumerate}

We have the following obvious

\begin{lemma}
$(i_k,r_k), (k=0,1,\cdots )$ is a finite sequence. In particular, the sequence stops at
$(i_m, r_m)$.
\end{lemma}

{\bf Proof.} By the construction of $i_k$'s, we have
\[
\sigma_{i_0}>\sigma_{i_1}>\cdots>\sigma_{i_k}>\cdots.
\]
Thus all $i_k$'s must be distinct. But $0\leq i_k\leq p$. 
So the length of the sequence is at most  $p+1$.

\qed

Let 
\begin{equation}\label{4-2}
\psi(x)=\underset{i\geq 0}{\rm Min} (\delta_i+\sigma_ix).
\end{equation}
The function $\psi(x)$ is a piecewise linear function, its derivative exists almost
everywhere. $r_k, (k=1,\cdots,m)$ are the non-smooth points of $\psi(x)$.

\begin{lemma}\label{lem42}
Assuming that $\sigma_{i_m}=0$, we have
\begin{equation}\label{4-3}
\sum_{k=0}^{m-1}(-\delta_{i_k}+\delta_{i_{k+1}})(
\sigma_{i_k}+\sigma_{i_{k+1}}-1)
=
\int_0^\infty\psi'(x)(\psi'(x)-1) dx.
\end{equation}
\end{lemma}

{\bf Proof.} First, let's remark that for $x$ large enough, $\psi\equiv\delta_{i_m}$
is a constant. Thus the integral in the lemma is convergent.

By definition of  $r_k (k=0,\cdots, m)$ in \eqref{4-1}, we have
\[
-\delta_{i_k}+\delta_{i_{k+1}}=(\sigma_{i_k}-\sigma_{i_{k+1}}) r_{k+1}
\]
for $k=0,\cdots, m-1$. Thus we have
\[
\sum_{k=0}^{m-1} (-\delta_{i_k}+\delta_{i_{k+1}})(\sigma_{i_k}
+\sigma_{i_{k+1}}-1)=\sum_{k=0}^{m-1}r_{k+1}
(\sigma_{i_k}^2-\sigma_{i_{k+1}}^2)+(\delta_{i_0}-\delta_{i_m}).
\]
The second term of the above equation is equal to 
\[
-\int_0^\infty\psi'(x) dx.
\]
For the first term, using the summation by parts, we have
\[
\sum_{k=0}^{m-1}r_{k+1}(\sigma_{i_k}^2-\sigma_{i_{k+1}}^2)=
r_1(\sigma_{i_0})^2+\sum_{k=1}^{m-1}\sigma_{i_k}^2(r_{k+1}-r_k)=
\int_0^\infty\psi'(x)^2 dx.
\]
Combining the above two equations, we get ~\eqref{4-3}.

\qed

Consider the smooth hypersurface $M\subset CP^n$ defined by the
polynomial $F=0$ of degree $d$. Let $X=\sum_{i=0}^n\lambda_iZ_i\frac{\pa}{\pa
Z_i}$ be the vector field for integers $(\lambda_0,\cdots,\lambda_n)$ such that
$\sum \lambda_i=0$. Let $M_t$  be defined by the equation
\begin{equation}\label{4-4-1-1}
F_t(Z_0,\cdots,Z_n)=F(t^{-\lambda_0}Z_0,\cdots, t^{-\lambda_n}Z_n).
\end{equation} 
We write $F_t$ as
\begin{equation}\label{4-4}
F_t=t^\delta\sum_{i=0}^pa_it^{\delta_i}Z_0^{\alpha^i_0}\cdots
Z_n^{\alpha_n^i},
\end{equation}
where $\delta_0=0$, and $\delta_i\geq 0, i=1,\cdots,p$. By~\eqref{4-4-1-1},
we have 
\begin{equation}\label{4-4-2}
X(Z_0^{\alpha_0^i}\cdots Z_n^{\alpha_n^i})
=-(\delta_i+\delta) Z_0^{\alpha_0^i}\cdots Z_n^{\alpha_n^i}
\end{equation}
for $i=0,\cdots, p$.

In what follows we assume that the choice of $(\lambda_0,\cdots,\lambda_n)$ is
``generic'' in the following sense:
\begin{enumerate}
\item All $\delta_i$'s are distinct;
\item None of the three lines defined by $\delta_i+\alpha_k^ix$ for $i=0,\cdots,p$
intersect at the same points, where $k=0,\cdots,p$.
\end{enumerate}
Without losing generality, we may assume that $\delta_0=0,
a_0=1$ and $0=\delta_0<\delta_1<\delta_2<\cdots<\delta_p$ .  We also assume that
$a_0,\cdots, a_p$ are all non-zero. Furthermore, since $M$ is smooth, we see
that for each $0\leq k\leq n$, there is an $0\leq i\leq p$ such that $\alpha_k^i=0$.

Let $U_i=\{[Z_0,\cdots,Z_n]\in CP^n|
|Z_i|>\frac 12|Z_j|, j=0,\cdots,n\}$. Then $\cup U_i=CP^n$. Let
$P_i=\{Z_i=0\}$ and $P_{ij}=P_i\cap P_j$ for $i\neq j$ and $i,j=0,\cdots, n$. Let
$\sigma>0$ be chosen so that $\sigma<\frac 1d\underset{i\geq 1}{\rm Min}
(\delta_i)$ (Note that $\underset{i\geq 1}{\rm Min}
(\delta_i)>0)$  and define
\[
V_{ij}^t=\{z|d(z,P_{ij})<|t|^\sigma\}, i\neq j, i,j=0,\cdots,n,
\]
where $d(\cdot,\cdot)$ is the distance induced by the Fubini-Study metric on $CP^n$.

By~\eqref{4-4}, we see that $t^{-\delta} F_t\rightarrow Z_0^{\alpha_0^0}\cdots
Z_n^{\alpha_n^0}$ as $t\rightarrow 0$. Intuitively, $M_t$ goes
to the hyperplanes defined by $Z_0^{\alpha_0^0}\cdots
Z_n^{\alpha_n^0}=0$. The following lemmas make
this observation rigid.

\begin{lemma}\label{lem33}
There is a $\sigma_1>\sigma$ such that for any  $0\leq k\leq n$ and
\[
[Z_0,\cdots,Z_n]\in (M_t-\cup_{i,j=0}^n V_{ij}^t)\cap U_k,
\]
one can find a unique $l\neq k$ such that

\[
\left|\frac{Z_l}{Z_k}\right|<|t|^{\sigma_1}
\]
for $t$ small enough.
\end{lemma}

{\bf Proof.} By~\eqref{4-4} we have
\begin{equation}\label{4-4-8}
|Z_0^{\alpha_0^0}\cdots Z_n^{\alpha_n^0}|
\leq
2^d\sum_{i=1}^p|a_i||t|^{\underset{i\geq 1}{\rm Min}(\delta_i)}|Z_k|^d.
\end{equation}
Thus if for any $l\neq k$,
\[
\left|\frac{Z_l}{Z_k}\right|\geq |t|^{\sigma_1},
\]
we could have
\[
|Z_0^{\alpha_0^0}\cdots Z_n^{\alpha_n^0}|\geq
|t|^{\sigma_1d}|Z_k|^d.
\]
This is a contradiction since we  choose $\sigma_1$ such that
\[
\sigma<\sigma_1<\frac 1d \underset{i\geq 1}{\rm Min}(\delta_i).
\]
 
\qed

We are now going to prove that
for $t$ small enough, the connected components of $M_t
\backslash\cup V_{ij}^t$ are graphs. We set 
\[
\tilde P_i=P_i-\cup_{j\neq i}V_{ij}^t,
\]
and let
\[
Q_i=\{[Z_0,\cdots,Z_n]| [Z_0,\cdots,Z_{i-1},\underset{i}{0},
Z_{i+1},\cdots, Z_n]\in \tilde P_i\},
\]
for $i=0,\cdots,n$. By~\eqref{pp-2} and ~\eqref{pp-3}, we have
\begin{equation}\label{p-1}
\phi(x_0,\cdots, x_n)=\underset{0\leq i\leq p}{\rm Min} (
\delta+\delta_i+\alpha_0^ix_0+\cdots
+\alpha_n^ix_n).
\end{equation}

\begin{rem}\label{normal}
$\phi$ and $\phi_i \,(i=0,\cdots,n)$ are defined even $\lambda_0,\cdots,\lambda_n$
are not choosing ``generically''. In the special case when
\[
XF=\kappa F,
\]
we have
\[
\phi_i(x)=-\kappa+(\underset{0\leq j\leq p}{\rm Min}\alpha_i^j)x
\]
for $0\leq i\leq n$. Thus if $M$ is a normal variety, we have
\[
\underset{0\leq j\leq p}{\rm Min}\alpha_i^j=0 \,\,{\rm or}\,\,1.
\]
In particular, in this case
\[
\phi_i'(x)(\phi_i'(x)-1)=0
\]
for $0\leq i\leq n$.
\end{rem}

\begin{prop}\label{prop1} Using  the notations as above, 
we have
\begin{align}\label{4-9}
\begin{split}
&\qquad\int_{M_t\cap Q_i}\sum_{A=0}^n\left(\frac{XF_t}{|\nabla F_t|^2}\right)_A
(\bar F_t)_A\omega^{n-1}\\&=
-\delta\alpha_i^{0}-\int_0^\infty \phi_i'(x) (\phi_i'(x)-1) dx+o(1),
\end{split}
\end{align}
for $i=0,\cdots, n$ as $t\rightarrow 0$.
\end{prop}

{\bf Proof.} 
For the sake of simplicity, we omit  unimportant constants in an inequality. Thus
in the proof of this proposition, $A\leq B$ means there is a constant 
$C$ independent of $t$ such that $A\leq CB$.

We just need to prove the theorem for the case $i=1$. 
If $\alpha_1^0=0$, then the proposition is automatically true since
$\phi_1'\equiv 0$. Thus
we assume that $\alpha_1^0\geq 1$.
We   work on $M_t\cap Q_1\cap U_0$ , without losing generality.

 We assume that
$(z_1,\cdots,z_n)=(\frac{Z_1}{Z_0},\cdots, \frac{Z_n}{Z_0})$ on $U_0$. Then
$F_t=0$ can be written as 
\begin{equation}\label{ff}
f=\sum_{i=0}^p a_i t^{\delta_i}z_1^{\alpha_1^i}\cdots z_n^{\alpha_n^i}=0
\end{equation}
with $a_0=1$ and $\delta_0=0$(see~\eqref{4-4}). The sequence $(\delta_i,\alpha_1^i),
(i=0,\cdots,p)$ is  assumed to be a 
``generic''
sequence mentioned at the beginning of this section.

For
$(z_1,\cdots,z_n)\in\tilde P_1\cap U_0$, we have
\[
|z_i|\geq |t|^\sigma,
\]
for $i=2,\cdots,n$.
Let $\xi_i^k (i=1,\cdots, \alpha_1^{i_k}-\alpha_1^{\alpha_{k+1}}, k=1,\cdots,m)$ be
the $(\alpha_1^{i_k}-\alpha_1^{\alpha_{k+1}})$-th roots of 
\[
-\frac{a_{i_{k+1}}}{a_{i_k}}t^{\delta_{i_{k+1}}-\delta_{i_k}}
z_2^{\alpha_2^{i_{k+1}}-\alpha_2^{i_k}}
\cdots
z_n^{\alpha_n^{i_{k+1}}-\alpha_n^{i_k}}.
\]
Then we have the following 

\begin{lemma} For $\sigma>0$ small enough,
there is a constant $\eps_0>0$ such that the solutions of $z_1$  of
$f=0$ satisfies
\[
|z_1-\xi_i^k|\leq|\xi_i^k|\cdot|t|^{\eps_0}
\]
for  $(i=1,\cdots, \alpha_1^{i_k}-\alpha_1^{i_{k+1}}, k=1,\cdots,m)$.
Furthermore, the balls $B_i^k=\{z\in \C||z-\xi_i^k|\leq|\xi_i^k||t|^{\eps_0}\}$
for $(i=1,\cdots, \alpha_1^{i_k}-\alpha_1^{\alpha_{k+1}}, k=1,\cdots,m)$ do not
intersect each other.
\end{lemma}

{\bf Proof.}  
In the proof, the scripts $i,k$ are always running in $(i=1,\cdots,
\alpha_1^{i_k}-\alpha_1^{\alpha_{k+1}}, k=1,\cdots,m)$, unless
otherwise stated.
We choose $\eps_1>0$ such that 
\[
\eps_1<\underset{0\leq k\leq m}{{\rm Min}}
\,
\underset{i\neq i_k,i_{k+1}}{{\rm Min}}
(\delta_i+\alpha_1^i r_k-\phi_1(r_k)).
\]
Define $f_k$ and $g_k$ as follows
\[
f_k=a_{i_k}t^{\delta_{i_k}}z_1^{\alpha_1^{i_k}}\cdots  z_n^{\alpha_n^{i_k}}    
+
 a_{i_{k+1}}t^{\delta_{i_{k+1}}}z_1^{\alpha_1^{i_{k+1}}}\cdots 
z_n^{\alpha_n^{i_{k+1}}},    
\]
and 
\[
g_k=f-f_k.
\]
By the definition of $\xi_i^k$, we have
\[
|t|^{r_k+C\sigma}\leq|\xi_i^k|\leq |t|^{r_k-C\sigma}
\]
for some constant $C$ independent of $t$.
We also have
\[
|g_k|\leq |t|^{\phi_1(r_k)+\eps_1-d\sigma}
\]
on $B_i^k$ and 
\[
|f_k|\geq |t|^{\phi_1(r_k)+\eps_0+d\sigma}
\]
on $\pa B_i^k$. We choose $\sigma$ small enough such that $\eps_1-d\sigma>\frac 34
\eps_1$ and $\eps_0$ small enough such that $\eps_0\leq \frac
14\eps_1$
Thus we have
\[
|f_k|>|g_k|
\]
on $\pa B_i^k$.
By the Rouch\'e Theorem, $f_k$ and $f=f_k+g_k$ have the same number of solutions
in
$B_i^k$. Since $f_k$ has only one solution in $B_i^k$, we prove
the first claim of the lemma. Next,
if $t$  is small enough, we have a constant $C$ such that
\[
|\xi_i^k-\xi_{i_1}^{k_1}|\geq C{\rm Max}(|\xi_i^k|,|\xi_{i_1}^{k_1}|).
\]
Thus if $t$ is small enough, $B_i^k$'s do not intersect each other.

\qed

{\bf  Continuation of the proof of Proposition~\ref{prop1}.} For simplicity, let
$F=F_t$. 
For fixed $i,k$, attaching the $B_i^k$ in the above lemma for each $p\in\tilde P_1\cap
U_0$,  we get a bundle $\tilde B_i^k$. On each $\tilde B_i^k$, by
~\eqref{4-4-2}, we have
\begin{align}\label{gg-1}
\begin{split}
&\qquad\qquad \sum_{A=0}^n\left(\frac{XF}{|\nabla F|^2}\right)_A
(\bar F)_A
=\frac{(XF)_1}{F_1}-\frac{(XF) F_{11}}{F_1^2}+o(1)\\
&\qquad=\frac{-(\delta+\delta_{i_k})\alpha_1^{i_k}+(\delta+\delta_{i_{k+1}})\alpha_1^{i_{k+1}}}
{\alpha_1^{i_k}-\alpha_1^{i_{k+1}}}\\&\qquad
-\frac{(-\delta_{i_k}+\delta_{i_{k+1}})
(\alpha_1^{i_k}(\alpha_1^{i_k}-1)-\alpha_1^{i_{k+1}}(\alpha_1^{i_{k+1}}-1))}
{(\alpha_1^{i_k}-\alpha_1^{i_{k+1}})^2}+o(1)\\
&=-\delta+
\frac{-\delta_{i_k}\alpha_1^{i_k}+\delta_{i_{k+1}}\alpha_1^{i_{k+1}}
+(\delta_{i_k}-\delta_{i_{k+1}})(\alpha_1^{i_k}+\alpha_1^{i_{k+1}}-1)}
{\alpha_1^{i_k}-\alpha_1^{i_{k+1}}}+o(1)
\end{split}
\end{align}
as $t\rightarrow 0$ for $k=0,\cdots,m-1$. 
By the same argument, the above equation is also true for $p\in\tilde P_1\cap U_l$
for $l\neq 0$. Thus the equation is true for $p\in\tilde P_1$. On the
other hand, by~\eqref{ff}, we have
\begin{equation}\label{gg-2} 
\det\pi=o(1)
\end{equation}
as $t\rightarrow 0$, where $\pi: Q_1\rightarrow\tilde P_1$ is the projection.
Thus by~\eqref{gg-1} and~\eqref{gg-2}, we have 
\begin{align*}
&\qquad\int_{M_t\cap Q_1}\sum_{A=0}^n\left(\frac{XF_t}{|\nabla F_t|^2}\right)_A
(\bar F_t)_A\omega^{n-1}\\&=
(-\delta\alpha_1^{0}+\delta_{i_m}
+\sum_{k=0}^{m-1}(\delta_{i_k}-\delta_{i_{k+1}})(\alpha_1^{i_k}
+\alpha_1^{i_{k+1}})){\rm vol} (CP^{n-1})+ o(1)
\end{align*}
as $t\rightarrow 0$, where $\alpha_1^m=0$ by the smoothness of
$M$. The proposition
follows from Lemma~\ref{lem42} and the fact ${\rm vol} (CP^{n-1})=1$.

\qed

\begin{lemma}\label{lem43}
Let $p$ be a fixed point in $M_t$ and let $d(x, p)$ be the distance from
$x\in CP^n$ to $x_0$ defined by  the Fubini-Study metric. Let $B_p(\eps)=\{
x\in CP^n|d(x,p)<\eps\}$. Then there is a constant $C$ independent of $p, \eps$
and $t$ such that
\begin{equation}\label{4-4-1}
\int_{M_t\cap B_p(\eps)}\omega^{n-1}
\leq
C\eps^{2n-2}(\eps^2\log |t|^{-1}+\log \eps^{-1}),
\end{equation}
for $t$ and $\eps$ small enough.
\end{lemma}

{\bf Proof.}  Consider the cut-off function
$\rho:\R\rightarrow
\R$ such that
$\rho\geq 0$ is smooth, $\rho\equiv 1$  on $[0,1]$ and $\rho\equiv 0$ on
$(-\infty, -1]\cup [2,+\infty)$. Then we have
\[
\int_{M_t\cap B_p(\eps)}\omega^{n-1}
\leq\int_{M_t}\rho(\frac{d(x,p)}{\eps})\omega^{n-1}.
\]
Let $F_t$ be the defining function of $M_t$. Then in the sense of 
distribution, we
have
\[
\bb\pa\bar\pa\log\frac{|F_t|^2}{(\sum_{i=0}^n|Z_i|^2)^d}=[M_t]-d\omega.
\]
Thus we have
\begin{equation}\label{4-5}
\int_{M_t}\rho\omega^{n-1}=d\int_{CP^n}\rho\omega^n+
\int_{CP^n}\rho\bb\pa\bar\pa\log
\frac{|F_t|^2}{(\sum_{i=0}^n|Z_i|^2)^d}\omega^{n-1}.
\end{equation}
We have an easy estimate for the first term of the right hand side of
~\eqref{4-5}:
\begin{equation}\label{4-6}
\int_{CP^n}\rho\omega^n\leq C\eps^{2n}.
\end{equation}
For the second term, assume that 
$p\in U_0=\{[Z_0,\cdots,Z_n]||Z_0|>\frac 12|Z_j|, j=1,\cdots,n\}$.  Then by
~\eqref{4-4}
\[
F_t=t^\delta Z_0^df_t , 
\]
where $f_t\rightarrow f_0=z_1^{\alpha_1^0}\cdots z_n^{\alpha_n^0}\not\equiv 0$.
Thus using integration by parts, we have
\begin{align}\label{4-7}
\begin{split}
&\int_{CP^n}\rho\bb\pa\bar\pa\log \frac{|F_t|^2}{(\sum_{i=0}^n|Z_i|^2)^d}
\omega^{n-1}\\&
\leq
C\eps^{2n}\log |t|^{-1}+\frac{C}{\eps^2}\left|\int_{|z|\leq 2\eps}\log |f_t| dV_0
\right|,
\end{split}
\end{align}
where $dV_0=(\bb)^n dz_1\wedge d\bar z_1\wedge\cdots\wedge d\bar z_n$ is the
Euclidean measure and $|z|^2=|z_1|^2+\cdots+|z_n|^2$. By changing the variables, the
second term of the above integral becomes
\begin{equation}\label{4-8}
\frac{C}{\eps^2}\int_{|z|\leq 2\eps}\log |f_t| dV_0=C\eps^{2n-2}\log \eps^{-1}
+C\eps^{2n-2}\left|\int_{|z|\leq 2}\log |f_t| dV_0\right|.
\end{equation}
By a theorem of Phong and Sturm~\cite{PS}, we have
\begin{equation}\label{4-9-101}
\int_{|z|\leq 2}\log |f_t|^{-1} dV_0 
\leq C
\end{equation}    
for $t$ small enough. 
~\eqref{4-4-1} follows from~\eqref{4-5}, ~\eqref{4-6},~\eqref{4-7},
~\eqref{4-8} and~\eqref{4-9}.

\qed

\begin{lemma}\label{lem46}
There exists a constant $C>0$ such that for $t$ small
\[
\sum_{i\neq j}\int_{V_{ij}^t\cap M_t}\omega^{n-1}\leq
C|t|^{2\sigma}\log |t|^{-1}.
\]
\end{lemma}

{\bf Proof.} Fixing $i,j$,
there is a constant $C_0$ independent of $\eps$ such that one can find points
$p_1,\cdots,p_m\in P_{ij}$  for $m=[\frac{C_0}{\eps^{2n-4}}]$, we have
\[
\overset{m}{\underset{k=1}{\cup}}B_{p_k}(\eps)\supset P_{ij}.
\]
Thus by the above lemma, we have
\[
\int_{V_{ij}^t\cap M_t}\omega^{n-1}\leq\sum_{k=1}^m\int_{M_t\cap B_{p_k}(
|t|^\sigma+\eps)}\omega^{n-1}.
\]
By the Lemma~\ref{lem43}, we have
\[
\int_{V_{ij}^t\cap M_t}\omega^{n-1}
\leq 
\frac{C}{\eps^{2n-4}}(|t|^\sigma+\eps)^{2n-2}
((|t|^\sigma+\eps)^2\log |t|^{-1}+\log (|t|^\sigma+\eps)^{-1}).
\]
The lemma follows from setting $\eps=|t|^\sigma$.

\qed

\begin{lemma}\label{lem48}
There exists a constant $C$ independent of $t$ such that
for any measurable subset $E$ of $M_t$
\[
\left|\int_E\pa\xi\wedge\bar\pa\theta\wedge\omega^{n-2}\right|\leq
C\sqrt{\log |t|^{-1}}\cdot\sqrt{{\rm meas}(E)},
\]
where the functions $\xi$ and $\theta$ are defined in~\eqref{3-1} and~\eqref{3-1-1},
respectively.
\end{lemma}

{\bf Proof.} 
Since $M_t$ is a submanifold, the Ricci curvature has an upper bound. Thus 
from ~\eqref{3-1}, we
have a constant $C$ such that
\begin{equation}\label{4-20} 
-\bb\pa\bar\pa\xi\leq C\omega
\end{equation}
On the other hand, since $[t^{\lambda_0}Z_0,\cdots,t^{\lambda_n}Z_n]\in M_t$ 
iff $[Z_0,\cdots,Z_n]\in M$, we have
\[
|\nabla F_t|^2(t^{\lambda_0}Z_0,\cdots,t^{\lambda_n}Z_n)
=\sum_{l=0}^n|t|^{-2\lambda_l}|F_k|^2(Z_0,\cdots,Z_n).
\]
Since $M$ is smooth, we have
\[
-C\log |t|^{-1}\leq |\xi|\leq C\log |t|^{-1}
\]
for some constant $C$.
Using integration by parts, from ~\eqref{4-20},  and the above estimate, we have
\[
\int_{M_t}|\nabla\xi|^2\omega^{n-1}\leq C\int_{M_t}(|\xi|+C\log |t|^{-1})
\omega^{n-1}\leq
C\log |t|^{-1}.
\]
If $E$ is a measurable subset of $M_t$, then we have
\[
\left|\int_E\pa\xi\wedge\bar\pa\theta\wedge\omega^{n-2}\right|\leq\int_E
|\pa\xi|\leq C\log |t|^{-1}\sqrt{{\rm meas}(E)},
\]
by the Cauchy inequality.

\qed

{\bf Proof of  Theorem~\ref{main}.} By Proposition ~\ref{prop1}, we have
\begin{align}\label{ll-1}
\begin{split}
&\quad\sum_{i=0}^n\int_{M_t\cap\cup Q_i}
\sum_{A=0}^n\left(\frac{XF_t}{|\nabla F_t|^2}\right)_A
(\bar F_t)_A\omega^{n-1}\\
&=-\delta d-\sum_{i=0}^n\int_0^\infty\phi_i'(x)(\phi_i'(x)-1)dx
+o(1)
\end{split}
\end{align}
as $t\rightarrow 0$.
We are going to prove that
\begin{equation}\label{4-22}
\int_{M_t\backslash
\overset{n}{\underset{i=0}{\cup}} Q_i}\sum_{A=0}^n\left(\frac{XF_t}{|\nabla
F_t|^2}\right)_A (\bar F_t)_A\omega^{n-1}=o(1)
\end{equation}
as $t\rightarrow 0$. In order to see this, let's recall that we have
\begin{align*}
&\qquad\int_{M_t\backslash
\overset{n}{\underset{i=0}{\cup}} Q_i}
\bb\pa\xi\wedge\bar\pa\theta\wedge\omega^{n-2}\\&
=\frac{1}{n-1}\left(-\int_{M_t\backslash
\overset{n}{\underset{i=0}{\cup}} Q_i}
\sum_{A=0}^n\left(\frac{XF_t}{|\nabla F_t|^2}\right)_A
(\bar F_t)_A\right.\\&+\left.\frac{\sum_{i=0}^n\lambda_i|(F_t)_i|^2}{|\nabla F_t|^2}
-(d-1)\theta\right)\omega^{n-1}
\end{align*}
by Lemma~\ref{lem310}.
Since $\theta$ and the function $\frac{\sum_{i=0}^n\lambda_i|F_i|^2}{|\nabla F|^2}$
are bounded, we have
\[
\left|\sum_{A=0}^n\left(\frac{XF_t}{|\nabla F_t|^2}\right)_A
(\bar F_t)_A\right|
\leq
\int_{M_t\backslash
\overset{n}{\underset{i=0}{\cup}} Q_i}
(|\pa\xi|+1)\omega^{n-1},
\]
by~\eqref{3-15-4}.
By Lemma~\ref{lem48} the righthanded side of the above equation is less than
or equal to
\[
C\sqrt{\log |t|^{-1}}\sqrt{{\rm meas} ({M_t\backslash
\overset{n}{\underset{i=0}{\cup}} Q_i})}+
{\rm meas} ({M_t\backslash
\overset{n}{\underset{i=0}{\cup}} Q_i}).
\]
If we can prove that there is a constant $C$ such that

\begin{equation}\label{4-21}
M_t\backslash
\overset{n}{\underset{i=0}{\cup}} Q_i
\subset
\underset{i\neq j}{\cup} V_{ij}^{Ct}.
\end{equation}
Then~\eqref{4-22}  will follow from Lemma~\ref{lem46}. To see ~\eqref{4-21}, let's
consider a point
$p\in{M_t\backslash
\overset{n}{\underset{i=0}{\cup}} Q_i}$. Without losing generality, we assume
that $p\in U_0$.
By ~\eqref{4-4-8}, we can find a $k\neq 0$ such that
\[
|Z_k|\leq |t|^\sigma|Z_0|
\]
for $t$ small enough. By definition, $p\notin Q_k$. Thus there is a $j\neq 0,k$ such
that
\[
|Z_j|\leq |t|^\sigma |Z_0|
\]
Thus $p\in V_{jk}^{Ct}$ for some constant $C$. ~\eqref{4-21} is proved.

Combining ~\eqref{ll-1} and ~\eqref{4-22}, we have
\[
\int_{M_t}
\sum_{A=0}^n\left(\frac{XF_t}{|\nabla F_t|^2}\right)_A
(\bar F_t)_A\omega^{n-1}=-\delta d-
\sum_{i=0}^n\int_0^\infty\phi_i'(x)(\phi_i'(x)-1) dx
+
o(1)
\]
as $t \rightarrow 0$.
Finally, 
since $\theta$ is a bounded function
\[
\int_{M_t}\theta\omega^{n-1}=\int_{M_0}\theta\omega^{n-1}+o(1)
\]
as $t\rightarrow 0$, where $M_0$ is defined as the zero set of
$Z_0^{\alpha_0^0}\cdots Z_n^{\alpha_n^0}=0$  counting  the multiplicity. In 
~\cite[Theorem 5.1]{Lu9},
it is proved that
\[
\int_{M_0}\theta\omega^{n-1}=-\frac{\delta}{n}.
\]
By~\eqref{m-2}, we have
\begin{align*}
&t\M'(t)=\frac{2}{d}\left(
\frac{\delta(n+1)(d-1)}{n}\right.+\left.\sum_{i=0}^n\int_0^\infty
\phi_i'(x)(\phi_i'(x)-1)dx\right)+o(1)
\end{align*}
as $t\rightarrow 0$ and Theorem~\ref{main} is proved.

\qed

\bibliographystyle{abbrv}
\bibliography{bib}

\end{document}